\documentclass[11pt,a4paper,final, leqno, twoside]{article}
\usepackage{amsmath}

\usepackage{amsthm}
\usepackage{amssymb}
\usepackage{amscd}
\usepackage{graphicx}
\usepackage{afterpage}

\begin{document}

\title    {\ \textbf{On invariant sets of diffeomorphisms} }
\author  {\small   Mehrzad Monzavi and Reza Mirzaei}
\date{}
\maketitle {\footnotesize
 \vskip
0.2 true cm \begin{abstract}
 We give a simple upper bound for the upper box
dimension of a backward invariant set of a $C^{1}$-diffeomorphism
of a Riemannian manifold. We also estimate an upper bound for the
box dimension of a forward invariant set of a $C^{1}$-mapping
with finite Brouwer degree in a Riemannian manifold.\\
{\it AMS Subject Classification: 58A05, 28A80.}

 \vskip 0.2
true cm \vskip 0.2 true cm {\it Key words:} Riemannian manifold,
Invariant set, Box dimension.
\end{abstract}}
 \vskip
0.2 true cm

\begin{flushleft}
    1. {\bf \small Introduction}
\end{flushleft}

The direct computation of the Hausdorff dimension of invariant
compact sets is a problem of high complexity. Therefore, it is
interesting to obtain analytic estimates of this dimension.
Recently, Many research studies have
 been developed on the investigation of the Hausdorff dimension of invariant
  sets of discrete dynamical systems.
  First results in this direction are given in [3] for compact subsets of
   $\mathbb{R}^{n}$ that
   are backward invariant under $C^{1}$-maps. Wolf in [24] gave Hausdorff
 dimension estimates
   (related to the values and behavior of $det D_{x}f$ and $\|D_{x}f\|$) for
   compact forward
   invariant sets of $C^{1}$-diffeomorphisms in $\mathbb{R}^{n}$. In [13, 14]
   the conditions of
   these estimates are weakened using a Lyapunov type function in $\mathbb{R}^{n}$. In [21],
    Temam gave upper bounds for
      the fractal dimension of flow-invariant sets in a Hilbert space, which is proved in [10] for
       vector fileds on
      Riemannian manifolds. In [7], Franz considered compact invariant sets of $C^{1}$-diffeomorphisms
      for which there
       exists an equivalent splitting of the tangent bundle. Qu and Zhou in [2] generalized the results
         of [24] to the map on smooth Rimannian manifolds with non-negative Ricci curvature. In [18],
         Wolf`s theorem is generalized
         to complete Riemannian manifolds without conditions on curvature. For
         further studies of
          estimation of upper bounds for  Hausdorff dimension of invariant compact sets,
          one may
           consult [4, 12, 15, 20, 23]. In the present paper, we estimate an upper bound
           for Hausdorff
            and box dimension of compact backward invariant sets of $C^{1}$-diffeomorphisms on Rimannian
            manifolds with conditions on $min S_{n}(D_{x}f)$ and $max |\det D_{x}f|.$ We also
             estimate an
             upper bound for the box dimension of a compact forward invariant set of a $C^{1}$-mapping
              with
              finite Brouwer`s degree on Riemannian manifolds.\\\\
               {\bf 2. Preliminaries}\\\\
               We will use the following definitions and facts.\\\\
(1) Consider a linear operator $L:E \rightarrow E'$ between two
Euclidean spaces of dimension $n$ with scalar products $<., .>_{E}$
and $<., .>_{E'}$ respectively. The adjoint operator of L is the
unique linear operator $L^{*}:E' \rightarrow E$ ,
 determined by $<Lx, y>_{E} =  <x, L^{*}y>_{E'}$ for all $x\in E$ and $y\in E'$.
  The eigenvalues of the positive semi-definite operator $\sqrt{L^{*}L}$ are the
  singular values of L. The singular values are all non-negative, usually listed in
   order to their size and multiplicity $ S_{1}(L)\geq S_{2}(L)\geq . . .\geq S_{n}(L)$ .
    The absolute value of the determinant of T is stated as the square root of the determinant
    of $T^{*}T$.\\\\
(2) Let M, N be  Rimannian manifolds of dimension $n$, $U$ be an open
subset of $M$.\\
(a) If $f:U\rightarrow M$  is
 a $C^{1}$-diffeomorphism.  We denote the tangent map of $f$
 at the point $x\in M$ by $D_{x}f: T_{x}M\rightarrow T_{f(x)}M$ and the norm of $f$ at that point, is defined
 by
 \[ \|D_{x}f\|=\sup\{|D_{x}f(v)| : v\in T_{x}M; |v|=1\}\]
(b) For $r > 0$ the $r$-neighborhood  of a set $F\subset M$ is
defined by
   \[B_{r}(F)=\{ x \in M:d(x,a)<\epsilon \ for \ some \ a\in A\}\]
(c) Let $F$ be a non-empty bounded subset of  $M$ and $N_{\delta}(F)$ be
the smallest number of balls of radius at
  most $\delta$ which can cover $F$. The upper box dimension of $F$ is defined(see [5]) by

\[ \overline{\dim}_{B}F=\limsup_{\delta\rightarrow
0}\frac{\log(N_{\delta}F)}{-\log\delta}\] (d) If $M$ and $N$ are
compact orientable  manifolds and $f: M \to N$ is a differentiable
map and $y\in N$ is a regular value of $f$, then the Brouwer
degree of $f$ at $y$ is defined (see [16]) by

\[ deg(f) = \sum_{x\in f^{-1}\in(y)} sgn(D_{x}(f))\]
Where,  $sgn(D_{x}(f))$ equals $+ 1$ or  $-1$ according to $D_{x}(f)$, which it preserves or reverses orientation.\\\\
(e) If $f:U \rightarrow M$ is a $C^{1}$-diffeomorphism onto its
image. A compact subset $K$ of $U$ is called
 forward $f$-invariant  if $f(K)\subset K$. If $K\subset f(K)$ then $K$ is called  backward $f$-invariant.\\\\
Authors of [1] gave a fractal dimension estimate for the invariant
set of a function $f:U \subset M \to M$ under the conditions
\[0<\min S_{n}(D_{x}f)<\sqrt{n}^{-1}\]          and
      \[(\max |det D_{x}f|) (\min S_{n}(D_{x}f))^{d-n} \leq 8^{-n}n^{\frac{-d}{2}}.\]
Qu and Zhou in [2] weakened the conditions and upgraded the
results of [1] under the condition that Ricci curvature is
non-negative. In this paper we generalize the results of [1] to
complete Riemannian manifolds (without considerations on Ricci
curvature). We will prove the following theorems.\\\\
{\textbf{Theorem 1.1.}} {\it Let $U\subset M$ be an open subset of a $n$-dimensional Riemannian
 manifold M, $f:U\rightarrow M$ be a $C^{1}$-diffeomorphism onto its image and $K$ be a backward $f$-invariant set.
 If $0<min S_{n}(D_{x}f)<1$  and
   there exists a number $d\in$ (0,d] such that
 \[(max |\det D_{x}f|) (min S_{n}(D_{x}f))^{d-n} \leq 1,\]
      then
      \[\overline{dim_{B}}K\leq d.\]
      }
The following theorem is proved in [24] under the condition $M=\mathbb{R}^{n}$. We prove  the same
result in more general case, when $M$ is  a complete Riemannian manifold.\\\\
{\textbf{Theorem 1.2.}} {\it Let $U\subset M$ be an open subset of a
 complete  Riemannian manifold $M$ and $f:U\rightarrow M$ be a $C^{1}$-mapping with the Brouwer degree $d $.
  Let  $K\subset U$ be a compact $f$-invariant set and suppose that $f$ has a non-zero Jacobian determinant. Put

\[ b=\lim_{m \rightarrow \infty} \frac{1}{m}\log(\min \{|\det D_{x}f^{m}|,
 x\in K\})\]
\[ s=\lim_{m \rightarrow \infty} \frac{1}{m}\log(\max \{\|D_{x}f^{m}\|
, x\in K\})\]
If $b>0$, then $s>0$ and
\[ \overline{\dim}_{B}K\leq n-\frac{b-\log d}{s}<n\]
\\\\}
{ \textbf{Remark 1.3.}} Because the Huasdorff dimension of a set is
 smaller or equal to its upper box dimension, theorems 1.1 and 1.2 also
 give upper bounds for the Huasdorff dimension of $K$.\\\\
{ \textbf{Remark 1.4.}} In Theorem 1.1, the assumption $0<min S_{n}(D_{x}f)<1$ would be
 unnecessary provided that the inequality $(max |\det D_{x}f|) (min S_{n}(D_{x}f))^{d-n} \leq 1$ was strict.\\\\
{\bf \large {2. Proofs of the theorems}}\\\\
\textbf{{Lemma 2.1}} (see [18]). {\it  If $K$ is a compact subset
of a Riemannian manifold $M$ and $dim M=n$, then
\[
\overline{dim}_{B}K \leq n+limsup_{r \rightarrow 0}
\frac{log(vol(B_{r}K ) )}{-log(r)}\]} \textbf{{Fact 2.2}} (see
[8, Theorem 2.92]). If  $M$ is a Riemannian manifold and
$x_{0}\in M$, then there is an open ball around $x_{0}$  such
that for any $x,y \in U$  there is a unique geodesic $\gamma$
joining $x$ to $y$ with the length equal to $d(x,y).$

\begin{flushleft}
\textbf{{Remark 2.3}} (see [2]). Let $B\subset U  $ be an open subset of a Riemannian
manifold $M$ and $f :U \to M$  a $C^{1}$-map. If $B$ is
bounded then \[vol(f(B))\geq inf_{x\in
B}|detD_{x}f|vol(B).\]
\end{flushleft}
\textbf{{Remark 2.4.}} If $U$ is an open subset of a Riemannian
manifold $M$ and  $f:U\rightarrow M$ a $C^{1}$-diffeomorphism on
its image, It is proved in [18] that if $K\subset U$ is a compact forward $f$-invariant set and
\[ b=lim_{m\rightarrow\infty} \frac{1}{m}log(min\{|detD_{x}f^{m}|,
 x\in K\})\]
\[ s=lim_{m\rightarrow \infty} \frac{1}{m}log(max\{\|D_{x}f^{m}\|
, x\in K\})\]
as well as $b>0$, then $s>0$ and
\[ \overline{dim}_{B}K\leq n-\frac{b}{s}<n\]
In a similar way we can prove the following theorem.\\\\
\textbf{{Theorem 2.5.}} {\it Let $U$ be an open subset of a Riemannian
manifold $M$ and  $f:U\rightarrow M$ a $C^{1}$-diffeomorphism on
its image. Let $K\subset U$ be a compact backward $f$-invariant set. Define
\[ b=lim_{m\rightarrow\infty} \frac{1}{m}log(min\{|detD_{x}f^{-m}|,
 x\in K\})\]
\[ s=lim_{m\rightarrow \infty} \frac{1}{m}log(max\{\|D_{x}f^{-m}\|
, x\in K\})\]
If $b>0$, then $s>0$ and
\[ \overline{dim}_{B}K\leq n-\frac{b}{s}<n\]
\\\\}
\textbf{Proof of Theorem 1.1.}\\\\
Since all norms in $\mathbb{R}^{n}$ are equivalent, the values of $b$ and $s$ are independent of the norm. Therefore the norm of $D_{x}f : T_{x}M\rightarrow T_{f(x)}M$ is equal to
        \[\|D_{x}f\| = \sqrt{\alpha_{n}}
        \]
 Where $\alpha_{1}$ is the maximum eigenvalue of $|D_{x}f|^{t}D_{x}f$.
Thus we have
  \begin{equation}
   S_{n}((D_{x}\varphi)^{-1}) = \|D_{x}\varphi\|^{-1}
  \end{equation}
    By (1) and the assumption that there exists a number $d\in$ (0,d] such that $(max |\det D_{x}f|) (min S_{n}(D_{x}f))^{d-n} \leq 1$, we have
  \begin{equation}
max|\det(D_{x}f)|\leq(min S_{n}(D_{x}f))^{n-d}= ((max\|(D_{x}f)^{-1}\|)^{-1})^{n-d}
  \end{equation}
   Using $f^{-1}(K)\subset K$ and (2) we have
     \[(max\|(D_{x}f)^{-1}\|)^{n-d} \leq min|det(D_{x}f^{-1})|\]
    Furthermore,
    \[min|det(D_{x}f^{-m})|=min|det(D_{x}f^{-1})\ldots det(D_{f^{-m+1}(x)}f^{-1})|\geq(min|det(D_{x}f^{-1})|)^{m}\]
    And
    \[max\|D_{x}f^{-m}\|\leq max(\|D_{x}f^{-1}\| \ldots \|D_{f^{-m+1}(x)}f^{-1}\|) \leq(max|det(D_{x}\varphi^{-1})|)^{m}\]
    By the assumptions of Theorem 1.1 we have
    \[max|det(D_{x}f)|\leq(min S_{n}(D_{x}f))^{n-d}\leq min S_{n}(D_{x})f\]
    Which results
    \[min det(D_{x}f^{-1})>1\]
    Therefore,

    \[\frac{b}{s}= lim_{m\rightarrow \infty}\frac{log(min{|detD_{x}\varphi^{-m}|:x\in k})}{log(max{\|D_{x}\varphi^{-m}\|:x\in k})} \geq \]\[
    \geq lim_{m\rightarrow \infty}\frac{log(min{|detD_{x}\varphi^{-1}|})^{m}}{log(max{\|D_{x}\varphi^{-1}\|})^{m}} \geq n-d \]
Now by Theorem 2.5 we have $\overline{dim}_{B}K <d.$ \begin{flushright}$\Box$ \end{flushright}
\textbf{Proof of Theorem 1.2.}\\\\
Since $f$ is a $C^{1}$-mapping, it follows from the definition of $b$ and $s$ and the continuity argument
that for each ${\delta}> 0$,  there exists $K_{\delta} \in \mathbb{N}$  and $\epsilon > 0$ such that
  \begin{equation}
    1<exp(k_{\delta}(b- \delta))< |detD_{x}f^{k_{\delta}}|
    \end{equation}
     and \begin{equation} \ \  ||D_{x}f^{k_{\delta}}||<
  exp(k_{\delta}(s+\delta))
  \end{equation}
for all $x \in B_{\epsilon}(K)$.
Since the Jacobian determinant is non-zero on $K$, there exists a neighborhood of $K$ on which $f$ can be considered locally as a $C^{1}$-diffeomorphism onto its image.
Let us assume $B_{\epsilon}(K)$ is such a neighborhood. It is possible to choose $\epsilon$  sufficiently small that for each $x \in K$ and
  each positive number $\varsigma \leq \epsilon $,
  $B_{\varsigma}(x)$ admits the results of Fact 2.2.
From now on consider the mapping $g=f^{k_{\delta}}$. Notice that $K$ is also backward $g$-invariant. Put
 \[r_{m}=\epsilon \ (expk(s+\delta))^{-m}< \epsilon\]
 and  \[B_{m} = B_{r_{m}}(K) \]

 for all $m\in \mathbb{N}.$\\\\
(1) Let $x\in B_{1}$, then there exists $y \in K$ such that $ d(x,y)< r_{_{1}} $. By Fact 2.2, there is a minimal geodesic $\gamma :[0,1] \to B_{1}(x)$  from $x$ to $y$. So
 \[ d(g(x),g(y))\leq
\int_{0}^{1}|\frac{d}{dt} g(\gamma(t))|dt \leq \int_{0}^{1}||D_{\gamma(t)}g||.|\frac{d}{dt} \gamma(t)|dt\leq\]\[
\leq (exp(k_{\delta}(s+\delta))) \int_{0}^{1}|\frac{d}{dt} \gamma(t)|dt = (exp(k_{\delta}(s+\delta))) d(x,y)<\]\[
< (exp(k_{\delta}(s+\delta))). r_{_{1}}=\epsilon\]\\\\
(2)  Now let $x\in B_{2}$, then there exists $y \in K$ such that $ d(x,y)< r_{_{2}} $. Similarly, we have
\[d(g^{2}(x),g^{2}(y))\leq \int_{0}^{1}|\frac{d}{dt} g^{2}(\gamma(t))|dt\leq \]\[\leq \int_{0}^{1}||D_{_{\gamma(t)}}g^{2}||.
|\frac{d}{dt} \gamma(t)|dt = \int_{0}^{1}||D_{_{\gamma(t)}}g||.||D_{_{g^(\gamma(t))}}g||.|\frac{d}{dt} \gamma(t)|dt \leq
\]\[\leq (exp(k_{\delta}(s+\delta)))^{2}\int_{0}^{1}|\frac{d}{dt} \gamma(t)|dt=(exp(k_{\delta}(s+\delta)))^{2} d(x,y)< \]\[
< (exp(k_{\delta}(s+\delta)))^{2}.r_{_{2}}=\epsilon\]
\\\\
By induction, for any $m \in \mathbb{N}$ and any $x \in B_{m}$, there exists $y \in K$ such that $d(g^{m}(x),g^{m}(y))\leq \epsilon$. Since $g^{m}(K) \subset K$ then $g^{m}(x) \in B_{\epsilon}K$. Thus $g^{m}(B_{m}(K)) \subset B_{\epsilon} K$ for all $m \in \mathbb{N}.$
Which resuls,
\begin{equation}
  vol (g^{m}(B_{m}(K))) \leq vol(B_{\epsilon}(K))
\end{equation}
$B_{m}(K)$ is bounded, so there exist balls $U_{k_{1}}, ...,U_{k_{n}}$ such that $g^{m}|_{U_{k_{1}}}$ is a $C^{1}$-diffeomorphism
onto its image and $B_{\epsilon}(K) \subset \bigcup_{i=1}^{n} U_{k_{l}}$.
Put
 \[V_{k_{1}} = U_{k_{1}} \cap B_{\epsilon}(K)\]
 \[V_{k_{i}} = (U_{k_{i}} \cap B_{\epsilon}(K)) \setminus \bigcup_{s=1}^{i-1} V_{k_{s}}\]
 By Remark 2.3 we have,
\begin{equation}
vol(g^{m}(V_{k_{i}})) \geq inf_{x\in V_{k_{i}}} |det D_{x}g^{m}| vol(V_{k_{i}})
\end{equation}
 By (3),   We have $exp(K_{\delta}(b-\delta)) < |det D_{x}g|$,   thus $exp(K_{\delta}(b-\delta))^{m} < |det D_{x}g^{m}|.$
 Therefore by (6), we have
\begin{equation}
 vol(V_{k_{i}}) \leq exp(K_{\delta}(b-\delta))^{-m} vol(g^{m}(V_{k_{i}}))
\end{equation}
The sets $V_{k_{i}}$ are pairwise disjoint and $\bigcup_{1}^{n}V_{k_{i}} = B_{\epsilon}(K)$. Using (5) and (7) we get
 \[vol(B_{\epsilon}(K)) = \sum_{i=1}^{n} vol(V_{k_{i}}) \leq \]
 \[\leq \sum_{l=1}^{n} exp(K_{\delta}(b-\delta))^{-m}vol(g^{m}(V_{k_{i}})) \leq\]
 \[\leq d^{km} exp(K_{\delta}(b-\delta))^{-m} vol(B_{\epsilon}(K)).\]
 Therefore,
\[limsup_{r \rightarrow
0}\frac{log(vol(B_{r}(K)))}{-log(r)}
=limsup_{r_{m}\rightarrow 0}\frac{log(vol
B_{m}(K))}{-log(r_{m})}\leq\] \[\leq lim_{m\rightarrow \infty}
\frac{log(d^{km} exp(K_{\delta}(b-\delta))^{-m}vol(B_{\epsilon}(K)))
}{-log(\frac{\epsilon}{exp(k_{\delta}(s+\delta))^{m}})}=-\frac{b-\delta-logd}{s+\delta}
\
\]
Since $\delta$ is arbitrary small,
then \begin{equation} \ \ limsup_{r \rightarrow
0}\frac{log(vol(B_{r}(K)))}{-log(r)}\leq -\frac{b-logd}{s} \
 \end{equation}
  Now by (8) and Lemma 2.1, we get
\[ \overline{dim}_{B}K\leq n-\frac{b-logd}{s}\] \begin{flushright}$\Box$ \end{flushright}
\vskip 0.4
true cm

  {\footnotesize

  \begin{center}

\end{center}

\bigskip
\begin{flushright}
\footnotesize {{\it
 Department of mathematics \ \ \ \ \ \ \ \\
 Faculty of sciences \ \ \ \ \ \ \ \ \  \ \ \ \ \\
 I. KH. International university(IKIU)\\
 Qazvin, Iran \ \ \ \ \ \ \ \ \ \ \ \ \ \ \ \ \\
  $mehrzad_{_{-}}$monzavi@yahoo.com \ \ \ \ \\
 $r_{_{-}}$mirzaioe@yahoo.com \ \ \ \ \ \ \ \ \ \\

    }}

\end{flushright}
}
\end{document}